\documentclass[11pt,cleanfoot]{asme2ej}

\usepackage{graphicx} 
\usepackage{epsfig}
\usepackage{fancyhdr}
\usepackage{setspace}
\usepackage{helvet}
\usepackage[hyphens]{url}

\usepackage{hyperref}   
\hypersetup{
	colorlinks=true,
	linkcolor=blue,
	citecolor=blue,
	urlcolor=blue,
}
\usepackage[square,numbers]{natbib}

\usepackage{amsmath}
\graphicspath{{./figure/}}

\title{Beam-like topologically interlocked structures with hierarchical interlocking}

\author{Ioannis Koureas
    \affiliation{Institute for Building Materials \\
        ETH Zurich \\
        Zurich, Switzerland
    }	
}

\author{Mohit Pundir
    \affiliation{Institute for Building Materials \\
        ETH Zurich \\
        Zurich, Switzerland
    }
}

\author{Shai Feldfogel
    \affiliation{Institute for Building Materials \\
        ETH Zurich \\
        Zurich, Switzerland
    }
}

\author{David S. Kammer
    \affiliation{Institute for Building Materials \\
        ETH Zurich \\
        Zurich, Switzerland \\
        e-mail: dkammer@ethz.ch
    }
}

\begin{document}

\maketitle    
\doublespacing
\begin{abstract}
{\it 
Topologically interlocked materials and structures, which are assemblies of unbonded interlocking building blocks, are promising concepts for versatile structural applications. They have been shown to exhibit exceptional mechanical properties, including outstanding combinations of stiffness, strength, and toughness, beyond those achievable with common engineering materials. Recent work has established a theoretical upper limit for the strength and toughness of beam-like topologically interlocked structures. However, this theoretical limit is only achievable for structures with unrealistically high friction coefficients; therefore, it remains unknown whether it is achievable in actual structures. Here, we demonstrate that a hierarchical approach for topological interlocking, inspired by biological systems, overcomes these limitations and provides a path toward optimized mechanical performance. We consider beam-like topologically interlocked structures that present a sinusoidal surface morphology with controllable amplitude and wavelength and examine the properties of the structures using numerical simulations. The results show that the presence of surface morphologies increases the effective frictional strength of the interfaces and, if well-designed, enables us to reach the theoretical limit of the structural carrying capacity with realistic friction coefficients. Furthermore, we observe that the contribution of the surface morphology to the effective friction coefficient of the interface is well described by a criterion combining the surface curvature and surface gradient. Our study demonstrates the ability to architecture the surface morphology in beam-like topological interlocked structures to significantly enhance its structural performance.

}
\end{abstract}



\section{Introduction}
\label{sec:introduction}

Topologically interlocked structures (TIS) are assemblies developed from specially-shaped building blocks that hold together without adhesives. The interfacial contact and the blocks' geometry provide the structural integrity of TIS. Since the blocks are not bonded, they can stick, slip and rotate with respect to each other. Biological materials like nacre use the concept of topological interlocking on the micro-scale allowing them to have better properties compared to conventional engineering structures. For instance, mechanical properties of TIS, such as toughness and strength, are highly adjustable and can be controlled by modifying a given structure's geometrical and material properties~\cite{Djumas2016, Djumas2017, Schaare2008, Mirkhalaf2018, Mirkhalaf2019, YazdaniSarvestani2019}.

TIS were first introduced by Dyskin et al.~\cite{dyskin2001toughening,dyskin2001new}, who presented the concept of using segmented tetrahedral blocks to form a structure. Apart from high toughness and strength, TIS enjoy many advantages such as fracture tolerance~\cite{Dyskin2003}, structural integrity even when a few blocks fail~\cite{Dyskin2003}, and customizability, which comes due to the extensive design space of the TIS~\cite{Estrin2021}. In terms of generalizing the geometric design space, TIS can be classified into three main categories; plate-like structures~\cite{Mirkhalaf2018, Mirkhalaf2019, Khandelwal2015}, beam-like structures~\cite{Dalaq2019, Dalaq2020, Odessa2022, Koureas2022} and curvilinear-like structures~\cite{Casapulla2019}. Here, we focus on beam-like structures to take advantage of the simple geometric configuration, which can facilitate their modeling and analysis using both computational and theoretical methods~\cite{Khandelwal2012, Khandelwal2014, Khandelwal2015}. This can be particularly useful for understanding the underlying physics of TIS and for developing predictive models of their behavior which can then be transferred to the more complex plate- and curvilinear-like TIS.

\begin{figure}[ht!]
    \begin{center}
    \includegraphics[width=15cm]{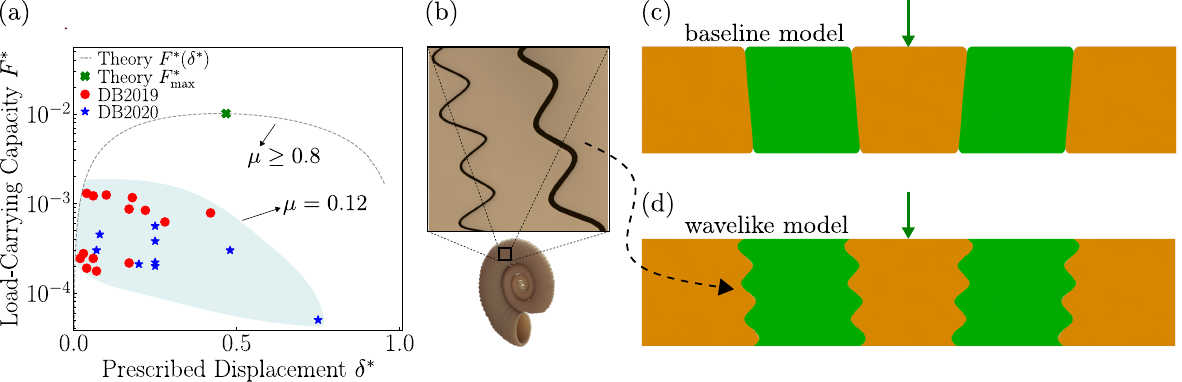}
    \caption{(a) Illustration of the normalized load-carrying capacity given by $F^* = F_y / (E \cdot h \cdot t)$ where $E$ is the Young's modulus, $h$ the height and $t$ the thickness of the structure, against the normalized prescribed displacement $\delta^* = \delta/h$. The dashed curve corresponds to the theoretical limit of planar-faced beam-like TIS with friction coefficient $\mu \geq 0.8$~\cite{Koureas2022}. Its maximum carrying capacity $F^*_\mathrm{max}$ (green cross) is compared to the maximum load-carrying capacity $F^*_{\mathrm{max}}$ of beam-like TIS from experiments by Dalaq and Barthelat \cite{Dalaq2019} (marked as DB2019) and Dalaq and Barthelat \cite{Dalaq2020} (DB2020). (b) Schematic of hierarchical suture interface in ammonite shell that inspired the current TIS design. (c-d) The rationale for the configurations: A beam-like topologically interlocked structure (TIS) with planar interfaces (c) as used in~\cite{Koureas2022}, and an updated configuration with hierarchical interfaces (d). The green arrows represent the loading direction.}
    \label{fig:timstructure}
    \end{center}
\end{figure}

The present work builds on a previous study~\cite{Koureas2022} in which we showed that beam-like TIS have a maximum theoretical load-carrying capacity. We demonstrated that this theoretical capacity is achieved by increasing the friction coefficient between the blocks up to the point where slip suppression is maximized. Similar observations that suggest forming the working envelop of TIS have been observed by Khandelwal et al.~\cite{Khandelwal2015} using fixed constraints and by Odessa and Shufrin~\cite{Odessa2022} through geometric modifications. In addition, an equivalent maximum theoretical limit for the load-carrying capacity has been observed for plate-like TIS~\cite{feldfogel2023scaling}. However, in all examined cases in~\cite{Koureas2022}, the friction coefficient that led to the maximum theoretical limit was much higher than the friction coefficients of common pairs of engineering materials. \autoref{fig:timstructure}a shows the theoretical limit for a planar-faced beam-like structure as shown in \autoref{fig:timstructure}c with friction coefficient $\mu \geq 0.8$ and compares it with other beam-like TIS from the literature~\cite{Dalaq2019, Dalaq2020}. We see at least one order of magnitude difference between the theoretical maximum load-carrying capacity obtained with $\mu \geq 0.8$ and the experimentally observed load-carrying capacities obtained with materials that have $\mu = 0.12$. This raises a question regarding the practical attainability of the theoretical maximum load-carrying capacity for beam-like TIS. Therefore, this study aims to show how the theoretical limit to the load-carrying capacity can be obtained with realistic friction coefficient values.

It has shown in~\cite{Li2011, Li2012, Lin2014, Malik2016, Malik2017, Malik2018, Khoshhesab2018, Wang2021} that hierarchical interlocking, of which surface waviness is one type, is a major contributor to the strength and the stiffness of biological materials like mollusk shells, conch shells, purpuratus shells, and ammonites shells (see \autoref{fig:timstructure}b). The concept of hierarchical interlocking by planar spikes over planar and non-planar interfaces was first introduced by Djumas et al. \cite{Djumas2017} in the context of plate-like TIS. Djumas et al. \cite{Djumas2017} demonstrated that the use of hierarchy provides an additional degree of interlocking and delays the onset of slip of the blocks. Inspired by these examples, we employ here a wave-based hierarchical interlocking to reach the theoretical capacities of TIS beams. Specifically, we introduce simple sinusoidal alterations to planar-faced surfaces (\autoref{fig:timstructure}c~and~d) to increase the frictional resistance between the blocks in beam-like structures. We explore the effects of these altered surfaces on the mechanical response for realistic values of friction coefficients and show that such simple geometric surface alterations are sufficient to increase the overall frictional resistance and thereby reach the maximum theoretical load-carrying capacity of the TIS, presented in Koureas et al. \cite{Koureas2022}.

In the following, we discuss the examined configuration and the surface alterations to increase the frictional resistance to reach the theoretical maximum load-carrying capacity in beam-like TIS (\autoref{sec:methodology}). In  \autoref{sec:results}, we present and discuss the effect of hierarchical morphology on the load-carrying capacity of TIS and give an outlook on how the gained knowledge can aid the design of TIS.

\section{Numerical Setup}
\label{sec:methodology}

\subsection{Examined Configuration}

\begin{figure}[ht!]
    \begin{center}
    \includegraphics[width=7.5cm]{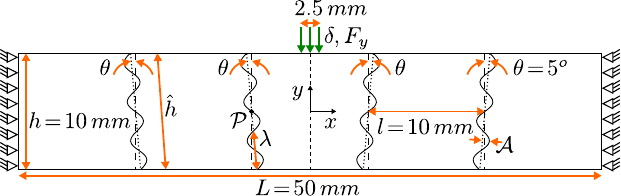}
    \caption{Schematic illustration of model set-up with the geometric parameters and boundary conditions in a beam-like TIS with hierarchical interfaces consisting of five blocks. Angle $\theta$ is computed between the vertical dashed lines and the inclined dotted lines.}
    \label{fig:timblocks}
    \end{center}
\end{figure}

To investigate the influence of surface morphology on the mechanical response of TIS, we focus on a centrally loaded beam-like structure with span length $L = 50~\textrm{mm}$, made of five blocks with height $h = 10~\textrm{mm}$ and length $l = 10~\textrm{mm}$ (see \autoref{fig:timblocks}). The central loading configuration was chosen because it constitutes the most critically loaded scenario. On top of the planar surfaces defined by $\theta = 5^o$, which were used by~\cite{Koureas2022} for obtaining the theoretical upper limit, we apply a hierarchical surface morphology. To generate surface morphology, we superimposed, a sinusoidal function defined by wavelength $\lambda$ and amplitude $\mathcal{A}$ on top of the inclined angle $\theta$.

The surface hierarchy of these beam-like TIS is modeled using a rotational (by angle $\theta$) transformation of a sinusoidal curve about the local midpoint $\mathcal{P}$ of the interface (see \autoref{fig:timblocks}) as
\begin{equation}
\label{eq:1}
    f(r) =
    \begin{bmatrix}
        \cos(\theta) & -\sin(\theta)\\
        \sin(\theta) & \cos(\theta)
    \end{bmatrix}
    \begin{bmatrix}
        r\\
        \mathcal{A} \sin\left(\left(r + \frac{h/2}{\cos(\theta)}\right) \frac{n\pi}{\frac{h/2}{\cos(\theta)}}\right)
    \end{bmatrix} +
    \begin{bmatrix}
        \mathcal{P}_x\\
        \mathcal{P}_y
    \end{bmatrix},
\end{equation}
where $n = \frac{\widehat{h}}{\lambda}$ is the number of oscillations across the interface, $\mathcal{A}$ is the amplitude, $\widehat{h}$ is the total length of the interface without the hierarchical modification (\autoref{fig:timblocks}), and $r$ is bounded by
\begin{equation}
\label{eq:bc}
    -\frac{h/2}{\cos(\theta)} \le r \le \frac{h/2}{cos(\theta)}~.
\end{equation}

For our study, the blocks are modeled as isotropic and linear elastic materials with Young's modulus $E$, Poisson's ratio $\nu$, and friction coefficient $\mu$. Further, we note that our model does not account for possible fracture of blocks.

\subsection{Parametric Study}
To analyze the effect of surface morphology, we design a range of hierarchical interfaces by varying $n$ and $\mathcal{A}$. The geometric parameters $h$ and $l$, as well as the material properties $E = 30~\mathrm{GPa}$ and $\nu = 0.2$, are kept constant.
Finally, $\mu$ ranges between $0.1-0.4$ covering a wide range of experimentally observed values for TIS~\cite{Djumas2017, Mirkhalaf2019, Dalaq2019, Dalaq2020, RezaeeJavan2020}.

To solve the proposed model, we use the Finite Element framework described in~\cite{Koureas2022}, which accounts for geometrical non-linearities~\cite{Bathe1975} and uses a penalty-based node-to-segment frictional contact algorithm~\cite{Konyukhov2004, Konyukhov2005, Laursen2002, Wriggers2007, Yastrebov2013, Zavarise2009, Zavarise2009a} to resolve contact and frictional constraints along interfaces. We employ Coulomb's friction law ($T = \mu N$, where $T$ and $N$ are tangential and normal interface tractions, respectively) to describe frictional forces along an interface and the stick-slip transition. Using the structures' symmetry about the $y$-axis, we model only the left half of the TIS (see \autoref{fig:timblocks}). The lateral surface of the TIS (\textit{i.e.}, left side) is constrained in $x$ and $y$ directions.

Finally, the load is applied as displacement-controlled boundary conditions along a section of the top boundary of the central block as described by
\begin{equation}
    \begin{aligned}
    &u_y(-\frac{l}{4} \le x \le 0, \frac{h}{2}) = -\delta~.
    \end{aligned}
\end{equation}
The 2D beam-like structures are modeled using first-order triangular elements under plane-strain conditions. The simulations are performed under static conditions. We note that we performed a mesh- and penalty-based sensitivity analysis to ensure that the load-displacement response and the interface behavior are converged.

\section{Results and Discussion}
\label{sec:results}

We explore the effect of $n$ and $\mathcal{A}$ on the global response using realistic values of $\mu$. Since we neglect fracture, we assume failure of the structure when the load drops to zero.

\subsection{Effect of surface morphology}

\begin{figure}[ht!]
    \begin{center}
    \includegraphics[width=15cm]{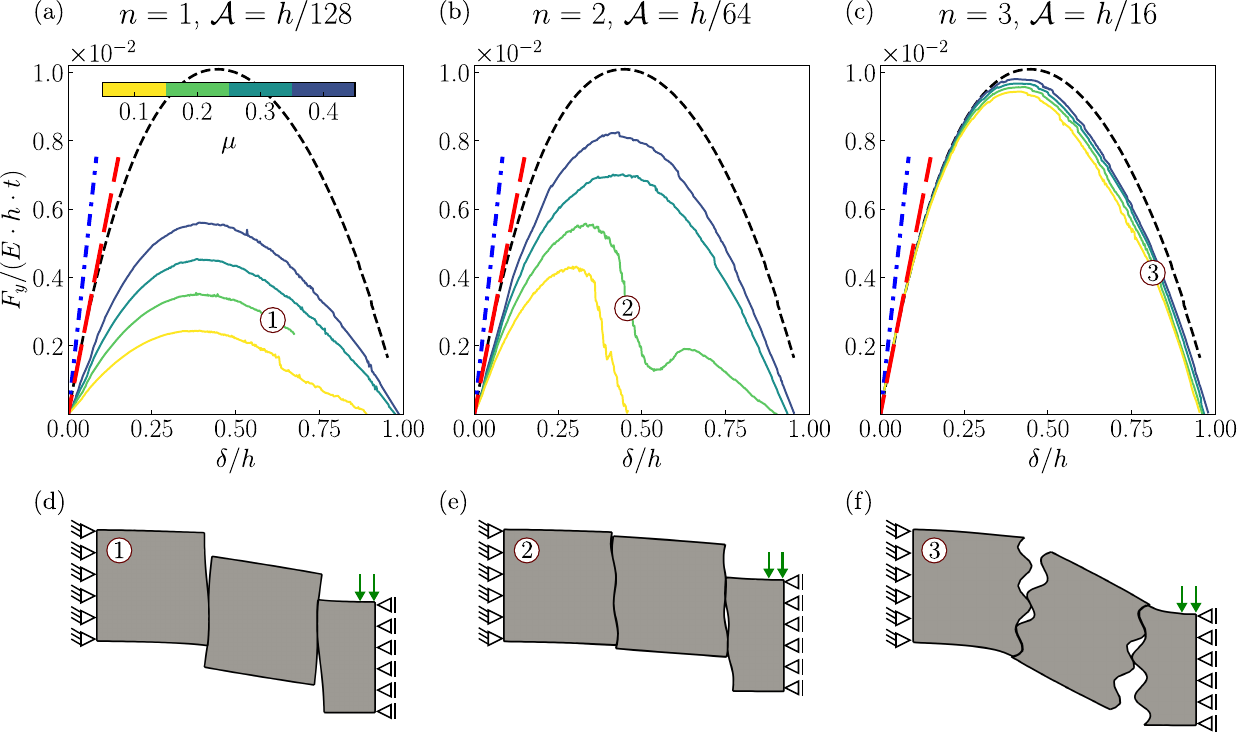}
    \caption{Normalized load-carrying capacity $F_y/(E \cdot h \cdot t)$ of the structure against the normalized prescribed displacement $\delta / h$ for various interface properties. (a-c) The curves correspond to structures with different friction coefficients $\mu$, oscillations $n$, and amplitude $\mathcal{A}$. (d-f) The schematics numbered 1, 2, and 3 show snapshots of the structures at given $\delta/h$. The curves are compared with the theoretical upper limit (dashed black curve) of planar-faced beam-like TIS with $\mu \geq 0.8$~\cite{Koureas2022}. The red dashed line represents the normalized structural stiffness $K$ of the saturated curve ($\approx 0.05$) and the blue dashed line represents $K$ of the monolithic equivalent (computed to be $\approx 0.089$ with finite-element simulations).}
    \label{fig:fdtrunc}
    \end{center}
\end{figure}

We analyze the normalized load-deflection curves $F_y - \delta$ for beam-like TIS with different friction coefficients $\mu$, number of oscillations $n$, and amplitude $\mathcal{A}$. To remove scalability effects, $F_y$ is normalized with respect to $h$, $t$, and $E$. $\delta$ is normalized with respect to $h$.
In all the considered cases, the load-carrying capacity increases as a function of $\mu$, $n$, and $\mathcal{A}$ (see \autoref{fig:fdtrunc}). The complete study for all values of $\mu$, $n$, and $\mathcal{A}$ is presented in \autoref{app:supplementary}. We can distinguish two different types of the $F_y - \delta$ curves. The first type (\autoref{fig:fdtrunc}a~and~c) follows a bell-shape behavior similar to other TIS~\cite{Mirkhalaf2019, Koureas2022, Dalaq2019}. In the second one (\autoref{fig:fdtrunc}b with $\mu = 0.2$), there are more than one local maxima. In all of our cases with more than one local maxima, the first maximum is always the global maximum.

\begin{figure}[ht!]
    \begin{center}
    \includegraphics[width=12cm]{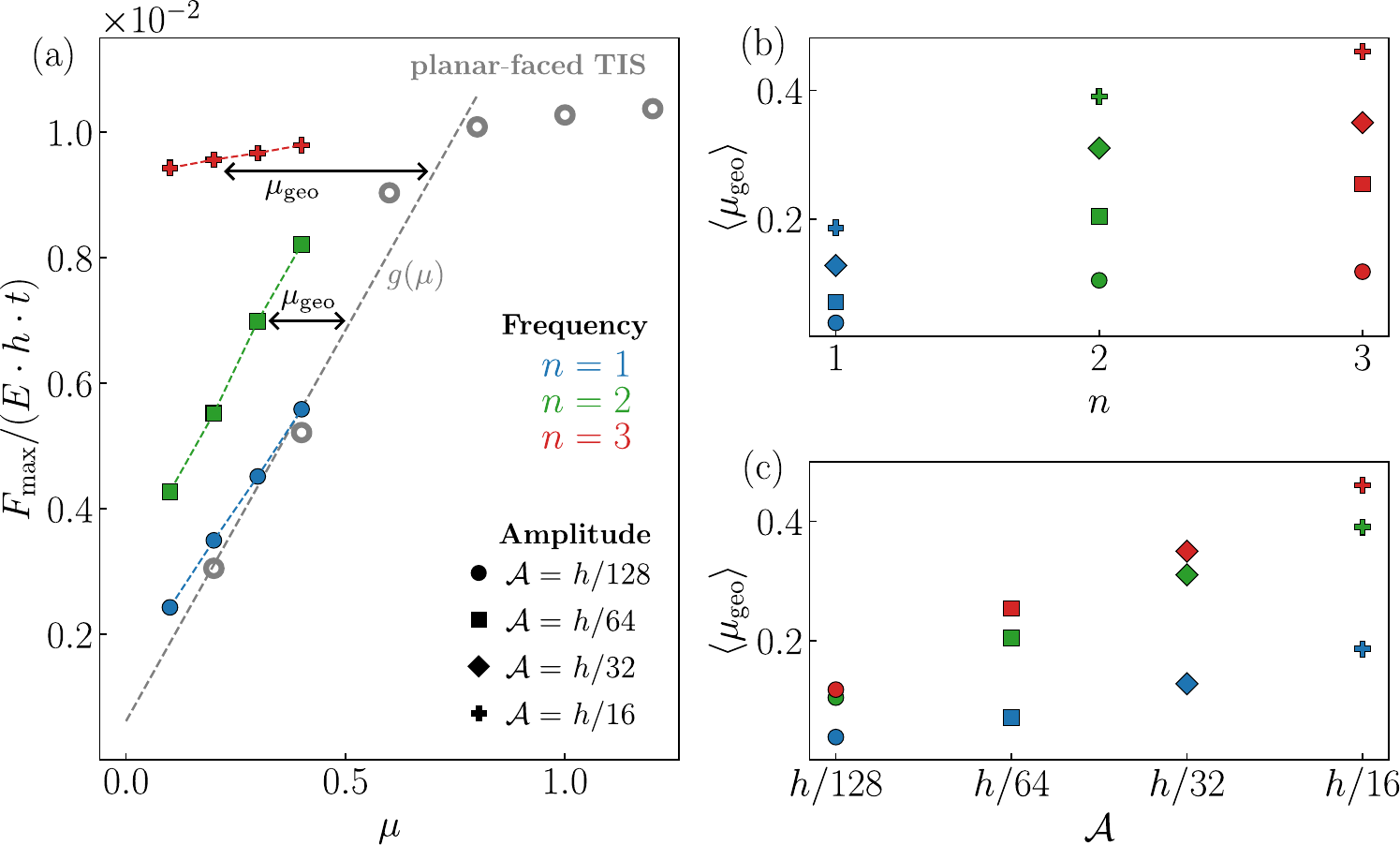}
    \caption{Geometric coefficient of friction. (a) Normalized maximum load-carrying capacity $F_\mathrm{max}/(E \cdot h \cdot t)$ against the friction coefficient $\mu$. The same cases as in \autoref{fig:fdtrunc} are shown, with oscillation $n = 1$ and amplitude $\mathcal{A} = h/128$ (shown in blue), $n = 2$ and $\mathcal{A} = h/64$ (shown in green) and $n = 3$ and $\mathcal{A} = h/32$ (shown in red). The grey points represent the normalized maximum load-carrying capacity for planar-face beam-like TIS, reported in~\cite{Koureas2022}. (b,c) Variation of average geometrical friction coefficient $\langle\mu_\mathrm{geo}\rangle$ with respect to $n$ and $\mathcal{A}$, respectively. Standard deviations are omitted here for clarity but shown in \autoref{fig:mu_geo_curv_rms}.}
    \label{fig:f_max_mu_amp_n}
    \end{center}
\end{figure}

The $F_y - \delta$ capacities are compared with the theoretical limit of the load-carrying capacity from our previous study with planar-faced beam-like TIS (dashed black line in \autoref{fig:fdtrunc}). We see that, by increasing $n$ and $\mathcal{A}$, the structures approach the theoretical limit with much smaller, and realistic friction coefficients. Moreover, for high values of $n$ and $\mathcal{A}$ (\textit{i.e.}, $n = 3$ and $\mathcal{A} = h/16$),  the effect of $\mu$ on the load-carrying capacity response becomes very small. This demonstrates that using hierarchical interlocking through architecting the surface morphology enables us to reach the theoretical capacities of TIS with realistic friction coefficients. In addition, \autoref{fig:f_max_mu_amp_n}a shows that the scaling of the structural load-carrying capacity with respect to the friction coefficient remains the same as observed in planar-faced TIS~\cite{Koureas2022,feldfogel2023scaling}. At the same level of normalized $F_\mathrm{max}$, the slope is the same between planar-faced and hierarchical TIS, but the level of friction (\textit{i.e.} $\mu$) at which this occurs is different. In the following, we will use this observation to characterize the effect of surface morphology on the mechanics of the interfaces.

Fundamentally, we note that the introduction of surface morphology increased the effective frictional resistance without increasing material friction, \textit{i.e.} for a given $\mu$ the normalized $F_\mathrm{max}$ is higher than in planar-faced TIS (see \autoref{fig:f_max_mu_amp_n}a). To describe this effect, we propose to distinguish between a material contribution $\mu$ and a geometrical contribution due to surface morphology $\mu_{\mathrm{geo}}$ to the effective friction resistance of the interfaces $\mu_{\mathrm{eff}}$ as experienced by the structure. Therefore, the effective friction coefficient can be defined as 
\begin{equation}
\mu_{\mathrm{eff}} = \mu + \mu_{\mathrm{geo}}~.
\label{eq:mu_eff}
\end{equation}
While in our previous work~\cite{Koureas2022} with planar-faced surfaces (\textit{i.e.}, $\mu_{\mathrm{geo}}=0$), the effective interface friction was given by $\mu_{\mathrm{eff}} = \mu$, and hence material friction had to be increased to unrealistically high values to reach the theoretical limit, for hierarchical interfaces, $\mu_{\mathrm{geo}} > 0$ provides new opportunities. In this case, $\mu$ remains in a range of realistic values for engineering materials, and by modifying the surface morphology of the structures, we can tailor $\mu_{\mathrm{eff}}$ via $\mu_{\mathrm{geo}}$. In other words, a geometric modification of the surface morphology can lead to the theoretical limit with realistic $\mu$ (see red points in \autoref{fig:f_max_mu_amp_n}a).

To get a better understanding of $\mu_{\mathrm{geo}}$ and its dependence on surface morphology (i.e., $n$ and $\mathcal{A}$) we compare the maximum load-carrying capacity $F_{\mathrm{max}}$ between the hierarchical TIS and the equivalent planar-faced beam-like TIS from Koureas \textit{et al.}~\cite{Koureas2022}. We describe $F_{\mathrm{max}}$ for planar-face beam-like TIS as a function of $\mu_{\mathrm{eff}}$, that is $F_{\mathrm{max}}^P = g(\mu_{\mathrm{eff}})$ (as represented by grey line in~\autoref{fig:f_max_mu_amp_n}, note that $\mu_{\mathrm{eff}} = \mu$ for planar-faced TIS). We then compute $\mu_{\mathrm{eff}}$ for hierarchical surfaces by projecting $F_{\mathrm{max}}$ for a given set of $(\mu,~n,~\mathcal{A})$ on to the function $g$, that is $\mu_{\mathrm{eff}} = g^{-1}(F_{\mathrm{max}}(\mu, n, \mathcal{A}))$. Using \autoref{eq:mu_eff}, we calculate the corresponding $\mu_{\mathrm{geo}}$ (as shown by arrows in \autoref{fig:f_max_mu_amp_n}a). We perform same calculations to compute $\mu_{\mathrm{geo}}$ for every combination of $n$, $\mathcal{A}$ and $\mu$ considered in the parametric study (see \autoref{app:supplementary}). An average geometrical friction coefficient $\langle\mu_\mathrm{geo}\rangle$ is then computed for a given set of values of $n$ and $\mathcal{A}$ (i.e., $\langle\mu_\mathrm{geo}\rangle = \sum\limits_{\mu_i}\mu_{\mathrm{geo}}(\mu_i, n,\mathcal{A}) /4$ where $\mu_i\in[0.1, 0.2, 0.3, 0.4]$). We observe $\langle\mu_\mathrm{geo}\rangle$ increases with increasing $n$ or $\mathcal{A}$ (see \autoref{fig:f_max_mu_amp_n}b and~c). However, none of these two parameters is sufficient to describe the observed $\mu_\mathrm{geo}$.

\begin{figure}[ht!]
    \begin{center}
    \includegraphics[width=\textwidth]{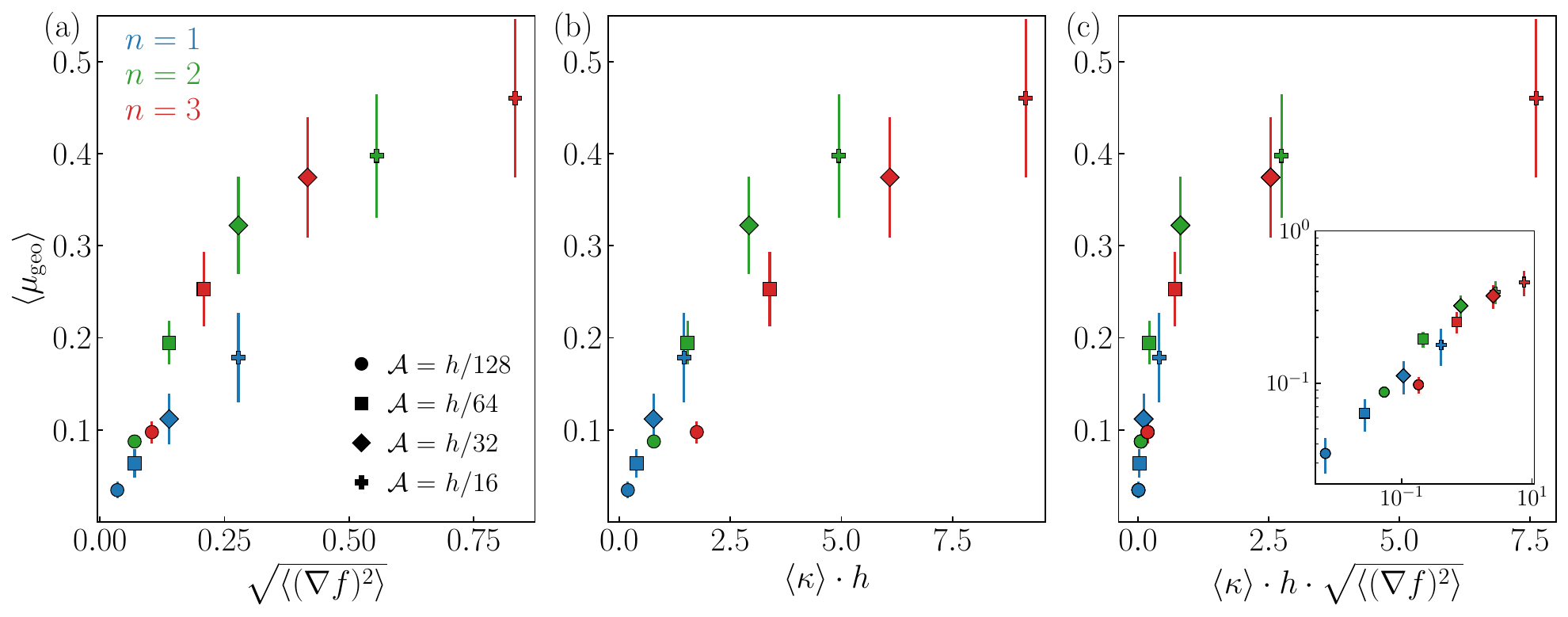}
    \caption{Variation of average geometrical friction coefficient $\langle\mu_\mathrm{geo}\rangle$ with respect to (a) root mean square of the surface gradient $\sqrt{\langle(\nabla f)^2\rangle}$, (b) normalized average curvature $\langle\kappa\rangle\cdot h$ and (c) product of $\langle\kappa\rangle\cdot h$ and $\sqrt{\langle(\nabla f)^2\rangle}$. Error bars indicate computed standard deviations. (inset) Same plot on log-log axis.}
    \label{fig:mu_geo_curv_rms}
    \end{center}
\end{figure}

The effects of $n$ and $\mathcal{A}$ can be combined in a measure that is computed via the gradient of the surface morphology, as commonly done in tribology~\cite{RanganathNayak1971, Yastrebov2012}. For this purpose, we compute the root mean square of the surface gradient $\sqrt{\langle(\nabla f)^2\rangle}$ and compare the variation of $\langle\mu_\mathrm{geo}\rangle$ as a function of $\sqrt{\langle(\nabla f)^2\rangle}$. We observe that the cases with $n=2,~3$ collapse onto a single curve (see \autoref{fig:mu_geo_curv_rms}a), however, the behavior of $\mu_\mathrm{geo}$ for $n=1$ is different. This suggests that the surface gradient is a good measure for morphologies that are tending to rougher profiles (which approach conditions common in tribological systems). Alternatively, in the context of TIS, surface morphology has previously been described by the surface curvature~\cite{Dalaq2020}. Here, we describe the local surface curvature as $\kappa = \frac{|\nabla^2 f|}{(1+\nabla f^2)^{3/2}}$ where $f$ represents the surface profile function (see \autoref{eq:1}). We then compare $\langle\mu_\mathrm{geo}\rangle$ with the normalized average curvature $\langle\kappa\rangle\cdot h$. The results show that $\langle\kappa\rangle\cdot h$ as a function of curvature collapses well for $n=1, 2$ but leads to a different dependence for $n=3$ (see \autoref{fig:mu_geo_curv_rms}b). This is consistent with observations by~\cite{Dalaq2020}, where curvature was sufficient as the considered morphologies were relatively smooth. 

Based on these observations, we propose a combined surface measure to characterize the geometrical friction coefficient of hierarchical interfaces, which is given by the product of surface curvature and surface gradient: $\langle\kappa\rangle\cdot h \cdot \sqrt{\langle(\nabla f)^2\rangle}$. The results, shown in \autoref{fig:mu_geo_curv_rms}c, illustrate a unique dependence of $\mu_{\mathrm{geo}}$ on this proposed surface measure. Therefore, the proposed surface measure appears to be able to account for the full spectrum of morphologies of hierarchical interfaces and can provide a better characterization of a surface's geometrical friction coefficient.

\subsection{Interface mechanisms in TIS with hierarchical interlocking}
\begin{figure}[ht!]
    \begin{center}
    \includegraphics[width=7.5cm]{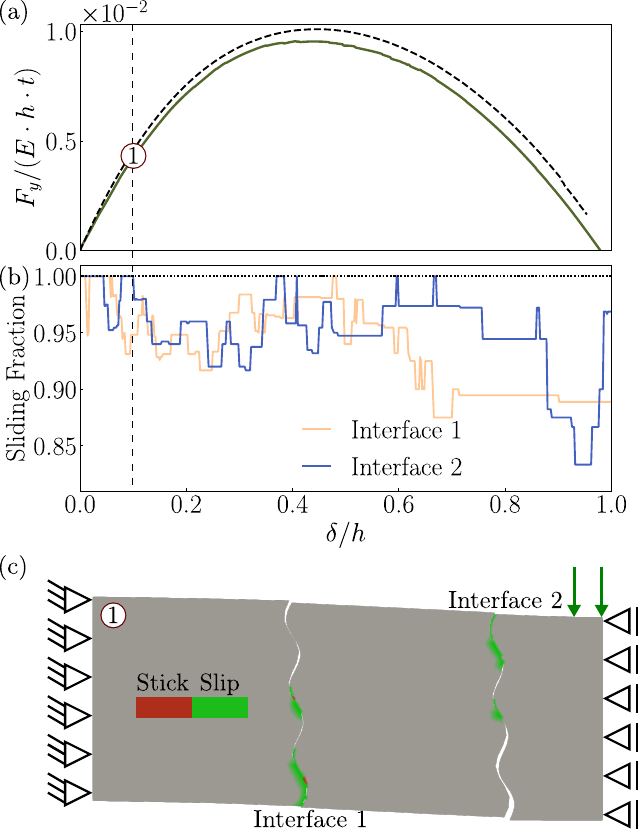}
    \caption{Interfacial behavior of representative case with friction coefficient $\mu = 0.4$, number of oscillations $n = 3$, and amplitude $\mathcal{A} = h/32$. (a) Normalized load-carrying capacity $F_y / (E \cdot h \cdot t)$. The curve is compared with the theoretical response (dashed curve) of planar-faced beam-like TIS with $\mu \geq 0.8$~\cite{Koureas2022}. (b) Sliding fraction against the prescribed displacement $\delta/h$ for each interface of the structure. (c) Snapshot capturing the sliding and sticking mechanisms at the interfaces of the structure at $\delta/h = 0.1$. The horizontal black dotted line instigates a fully sliding interface (sliding fraction $=$ 1).}
    \label{fig:stickslip}
    \end{center}
\end{figure}

The origin of the observed saturation level can be understood as the result of a collective interface behavior. Therefore, to analyze the impact of local interface mechanisms on structural properties, we define the sliding fraction as the contact area that slides at a given interface. \autoref{fig:stickslip} shows that the structures can approach the theoretical limit of the load-carrying capacity (\autoref{fig:stickslip}a) even when some interfaces are fully sliding (\textit{i.e.}, sliding fraction $\approx$ 1)~(\autoref{fig:stickslip}b~and~c), compared to~\cite{Koureas2022} where the saturation was reached only when all interfaces stick.

\begin{figure}[ht!]
    \begin{center}
    \includegraphics[width=7.5cm]{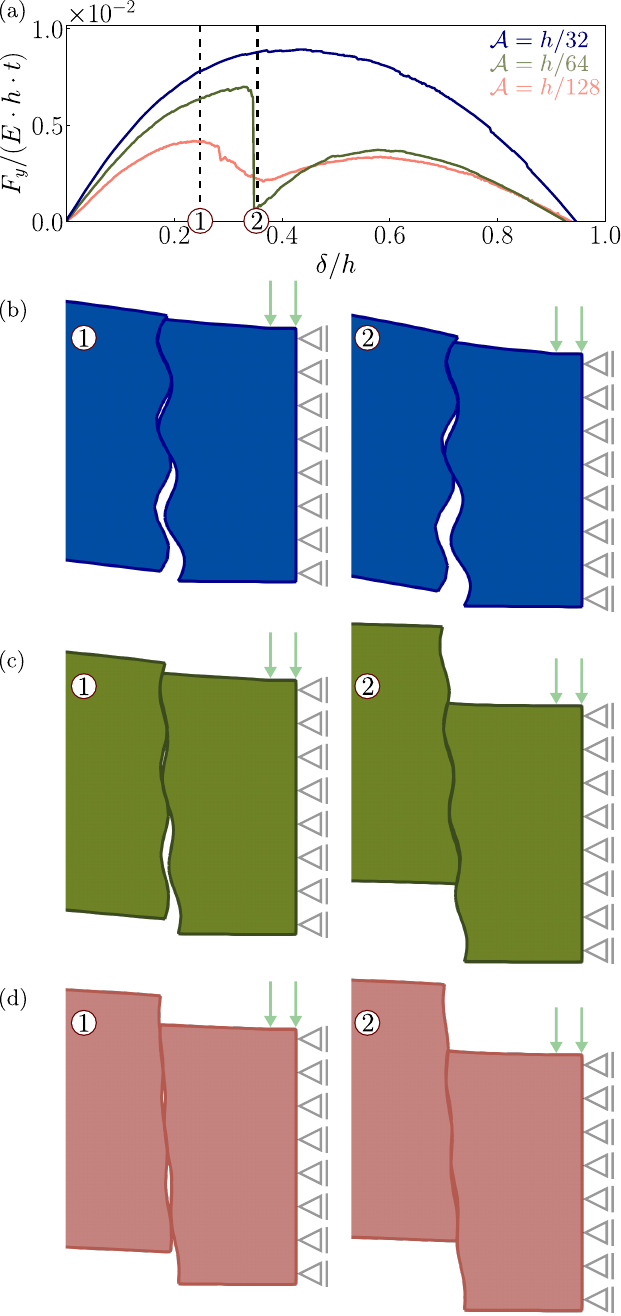}
    \caption{(a) Normalized load-carrying capacity $F_y/(E \cdot h \cdot t)$ against the normalized prescribed displacement $\delta/h$ of beam-like TIS with friction coefficient $\mu = 0.2$, number of oscillations $n = 3$ and amplitudes $\mathcal{A} = h/32,~h/64,~h/128$. (b-d) Snapshots at $\delta/h = 0.25$ and $\delta/h = 0.35$ show the evolution of the central block for (b) a case where the structure reaches the theoretical limit, (c) a case where interfacial sliding presents an instability marked by a sharp load drop, and (d) a case with smooth and continuous sliding of the central block.}
    \label{fig:truss_model_wavy}
    \end{center}
\end{figure}

In hierarchical TIS, local sliding dominates the interface mechanisms, which has direct implications on the structural behavior.
Our simulations show that we can distinguish three types of $F_y - \delta$ curves (see \autoref{fig:truss_model_wavy}). The first type (with $\mathcal{A} = h/32$), is able to reach the theoretical limit (see \autoref{fig:truss_model_wavy}a and \autoref{fig:truss_model_wavy}b) thanks to an effective combination of $\mu$ and $\mu_{\mathrm{geo}}$. The second type (with $\mathcal{A} = h/64$), presents a sharp load drop that indicates structural instability. From \autoref{fig:truss_model_wavy}a and \autoref{fig:truss_model_wavy}c, it is clear that this drop is associated with an interfacial slip-jump across a full wave. This unstable type of response corresponds to an undesirable combination of $\mu$ and $\mu_{\mathrm{geo}}$. In the third type (with $\mathcal{A} = h/128$), the structure never reaches the theoretical limit. The low value of $\mu_{\mathrm{eff}}$ makes the central block to smoothly slide across the interface (see \autoref{fig:truss_model_wavy}a  and \autoref{fig:truss_model_wavy}d) providing a stable movement of the block. These results suggest that a hierarchical modification of a beam-like TIS can potentially present these three regimes, which need to be identified for a given structure to avoid falling into the unstable regime which might cause uncontrolled response behavior.

\subsection{General Comments}
In this work, we considered TIS as beam-like structures. To put the observed saturated response in perspective to a monolithic structure, we compare the initial normalized stiffness $K$ of saturated TIS with that of a monolithic equivalent (see \autoref{fig:fdtrunc}). We define $K$ for TIS by the secant slope in the $F_y - \delta$ curves at $\delta/h = 0.05$. The initial stiffness of a monolithic structure represents an upper limit for the structural stiffness of TIS (see blue dashed line in \autoref{fig:fdtrunc}). We observe that the saturated initial stiffness of the beam-like TIS is approximately $56\%$ of the stiffness of the monolithic beam. We also note that the maximum carrying capacity and toughness of monolithic beams depend directly on the material failure criterion, which is not defined in this study. Hence, strength and toughness comparisons between TIS and monolithic beams are left for future work.

Here, we focused on \emph{centrally} loaded beam-like TIS as it is the most critical load configuration and, hence, most commonly studied in the literature. Other boundary conditions, where the load is asymmetrically applied on other blocks than the central one, will present the same mechanics but with quantitative differences. For instance, the initial stiffness and the maximum load-carrying capacity of the structure increase for asymmetric loading compared to the symmetric case, as can be expected from simple beam theory and the truss model by Khandelwal et al.~\cite{Khandelwal2012}. However, structures with interface morphology and friction that lead to the sticking of all interfaces will reach a configuration-specific saturation of the maximum load-carrying capacity. Hence, the findings of our work are qualitatively applicable to asymmetric loading configurations of TIS.

Another important structural property of TIS is the loading energy which is computed by the integral of the load versus imposed displacement (\textit{i.e.} the area under a $F_y - \delta$ curve). Maximizing the loading energy is a typical objective when designing TIS. For beam-like TIS, the maximum loading energy is directly related to the maximum load-carrying capacity, as has previously been demonstrated~\cite{Khandelwal2012, Khandelwal2014, Khandelwal2015}. Since the displacement $\delta$ cannot exceed the height $h$ in beam-like TIS, the saturated load-carrying capacity occurs for systems that reach $\delta = h$ (see \autoref{fig:fdtrunc}). Therefore, a saturation of the load-carrying capacity implies saturation of the loading energy for beam-like TIS.

Finally, it is important to note that our model does not consider the possibility of block fracture, which could potentially impact the performance of the TIS. To obtain the full theoretical range of the $F_y - \delta$ curves, we assumed that the failure of the structure occurs when the central block is pushed out of the assembly. While this scenario is realistic for many configurations~\cite{Malik2016, Malik2017, Malik2018, Wang2021}, it may not be applicable in all cases and leaves open questions regarding the exact point at which fracture initiation occurs, and how it will affect the structural response. Interestingly, local fracture in a block may not lead to structural failure, but it can introduce new smaller bodies that affect the overall response capacity.

\section{Conclusion}
\label{sec:conclusion}

This study focused on the effects of bio-inspired hierarchical interlocking on beam-like topologically interlocked structures (TIS). Specifically, we sought to test whether it was possible, by architecting the surface morphology of blocks, to reach the theoretical capacities of these structures with realistic friction coefficients, which cannot be achieved with planar-faced blocks. To this end, we tested a centrally loaded TIS, where different patterns of surface waviness were superimposed on the originally planar-faced blocks. We found that, with certain waviness patterns, it is possible to reach the theoretical limit at much lower and more realistic friction coefficients. Moreover, we found that the load-carrying capacity of structures with wavy block surfaces converges to that of a planar-faced beam-like TIS at high friction. Our observations demonstrate the usefulness of hierarchical interlocking in TIS and highlight its potential in the design of TIS.



%
% ...

\bibliographystyle{asmems4}

\bibliography{asme2e}

\newpage

\section{Supplementary Data}
\label{app:supplementary}

\begin{figure}[ht!]
    \begin{center}
    \includegraphics[width=14cm]{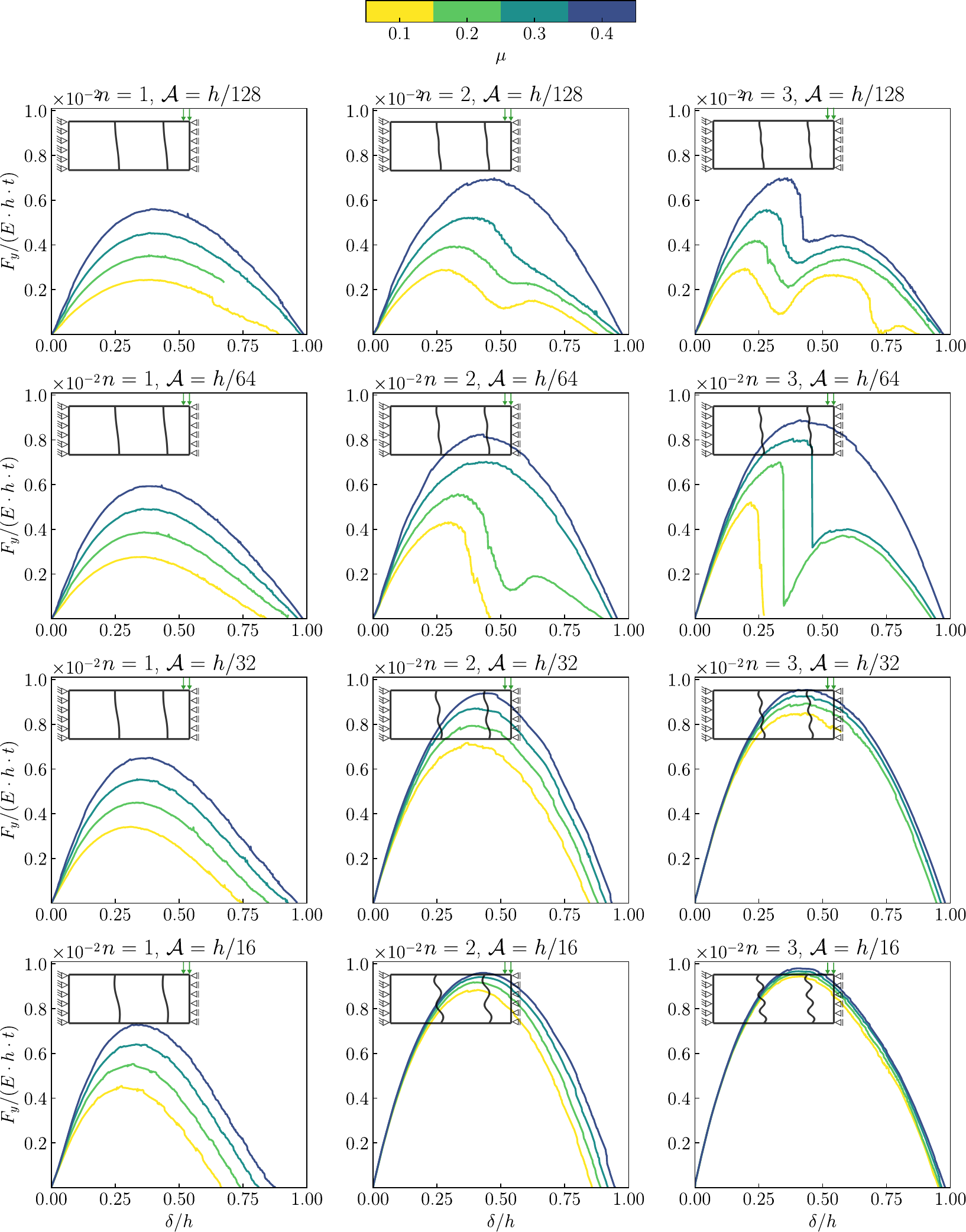}
    \caption{Normalized load-carrying capacity $F_y/(E \cdot h \cdot t)$ against the normalized prescribed displacement $\delta/h$. The curves correspond to structures with different friction coefficients $\mu$,  oscillations $n$, and amplitude $\mathcal{A}$. The symmetric model of the examined structure is shown in the top left corner of each sub-figure.}
    \label{fig:fdtrunc0}
    \end{center}
\end{figure}




\end{document}